\documentclass[12pt]{article}
\usepackage{amssymb,amsmath,latexsym}

\oddsidemargin=\evensidemargin
\advance\topmargin by-\headheight
\def\AGL{\mathop{\rm AGL}\nolimits}
\def\aut{\mathop{\rm Aut}\nolimits}

\def\End{\mathop{\rm End}\nolimits}
\def\F{\mathbb F}

\def\G{{\cal G}}
\def\GL{\mathop{\rm GL}\nolimits}
\def\GF{\mathop{\rm GF}\nolimits}
\def\GR{\mathop{\rm GR}\nolimits}

\def\id{\mathop{\rm id}\nolimits}

\def\ker{\mathop{\rm ker}\nolimits}
\def\lg{\langle}
\def\M{{\cal M}}
\def\Mat{\mathop{\rm Mat}\nolimits}

\def\N{{\mathbb N}}
\def\ov{\overline}
\def\O{{\cal O }}
\def\P{{\cal P }}
\def\PSL{\mathop{\rm PSL}\nolimits}
\def\R{{\cal R}}
\def\rad{\mathop{\rm Rad}\nolimits}
\def\rg{\rangle}

\def\SL{\mathop{\rm SL}\nolimits}

\def\sym{\mathop{\rm Sym}\nolimits}
\def\T{{\cal T}}

\def\ZZ{{\mathbb Z}}

\def\proof{{\bf Proof}.\ }
\def\bull{\vrule height .9ex width .8ex depth -.1ex }

\makeatletter
\renewcommand{\subsection}{\@startsection{subsection}{2}{0mm}{-2mm}{-2mm}
{\bf\normalsize}}

\def\sbsnt#1{\subsection{\hspace{-3mm}#1}}
\makeatother

\newtheorem{formula}{}[section]
\newtheorem{proposition}[formula]{Proposition}
\newtheorem{definition}[formula]{Definition}
\newtheorem{corollary}[formula]{Corollary}
\newtheorem{remark}[formula]{Remark}
\newtheorem{lemma}[formula]{Lemma}
\newtheorem{theorem}[formula]{Theorem}

\newtheorem{assumption}[formula]{Assumption}

\def\thrm{\begin{theorem}}
\def\thrml#1{\begin{theorem}\label{#1}}
\def\ethrm{\end{theorem}}
\def\rmrk{\begin{remark}}
\def\rmrkl#1{\begin{remark}\label{#1}}
\def\ermrk{\end{remark}}
\def\dfntn{\begin{definition}}
\def\dfntnl#1{\begin{definition}\label{#1}}
\def\edfntn{\end{definition}}
\def\nmrt{\begin{enumerate}}
\def\enmrt{\end{enumerate}}
\def\tm#1{\item[{\rm (#1)}]}
\def\qtn{\begin{equation}}
\def\qtnl#1{\begin{equation}\label{#1}}
\def\eqtn{\end{equation}}
\def\lmm{\begin{lemma}}
\def\lmml#1{\begin{lemma}\label{#1}}
\def\elmm{\end{lemma}}
\def\crllr{\begin{corollary}}
\def\crllrl#1{\begin{corollary}\label{#1}}
\def\ecrllr{\end{corollary}}
\def\ssmptn{\begin{assumption}}
\def\ssmptnl#1{\begin{assumption}\label{#1}}
\def\essmptn{\end{assumption}}
\def\prpstn{\begin{proposition}}
\def\prpstnl#1{\begin{proposition}\label{#1}}
\def\eprpstn{\end{proposition}}

\def\prblm#1#2{\vspace{2mm}{\bf\noindent #1.} {\it #2}\vspace{2mm}}

\begin{document}
\title{Constructions in public-key cryptography over matrix groups}
\author{
Dima Grigoriev \\[-1pt]
\small IRMAR, Universit\'e de Rennes \\[-3pt]
\small Beaulieu, 35042, Rennes, France\\[-3pt]
{\tt \small dima@math.univ-rennes1.fr}\\[-3pt]
\small http://name.math.univ-rennes1.fr/dimitri.grigoriev
\and
Ilia Ponomarenko
\thanks{Partially supported by RFFI, grants,  03-01-00349,
NSH-2251.2003.1. The paper was done during the stay of the author
at the Mathematical Institute of the University of Rennes.}\\[-1pt]
\small Petersburg Department of V.A.Steklov\\[-3pt]
\small Institute of Mathematics\\[-3pt]
\small Fontanka 27, St. Petersburg 191023, Russia\\[-3pt]
{\tt \small inp@pdmi.ras.ru}\\[-3pt]
\small http://www.pdmi.ras.ru/\~{}inp
}
\date{07.06.2005}
\maketitle

\begin{abstract}
The purpose of the paper is to give new key agreement protocols
(a multi-party extension of the protocol due to Anshel-Anshel-Goldfeld and a generalization
of the Diffie-Hellman protocol from abelian to solvable groups) and a new homomorphic public-key
cryptosystem. They rely on difficulty of the conjugacy and membership problems for subgroups
of a given group. To support these and other known cryptographic schemes we present a
general technique to produce a family of instances being matrix groups (over finite commutative
rings) which play a role  for these schemes similar to the groups $Z_n^*$ in the existing
cryptographic constructions like RSA or discrete logarithm.
\end{abstract}

\section*{Introduction}
One of the oldest cryptographical problems consists in constructing of a key agreement
protocol. Roughly speaking it is a multi-party algorithm, defined by a sequence of steps,
specifying the actions of two or more parties in order a shared secret becomes available to
two or more parties. Probably the first such procedure based on abelian groups is due to
Diffie-Hellman one (see~\cite{GB}). In fact, it concerns automorphisms of abelian (even
cyclic) groups induced by taking to a power. Some generalizations of this protocol to
non-abelian groups (in particular, the matrix groups over some rings) were suggested
in~\cite{PKHK} where security was based on an analog of the discrete logarithm problems in
groups of inner automorphisms. Certain variations of the Diffie-Hellman systems
over the braid groups were described in~\cite{KL}; there several trapdoor one-way
functions connected with the conjugacy and the taking root problems in the
braid groups were proposed. Recently, a general scheme for constructing key agreement
protocols based on algebraic structures was proposed in \cite{AAG}. In principle, it enables
us to construct such protocols for non-abelian groups and their automorphisms induced by
conjugations. In this paper we generalize to the non-abelian case the Diffie-Hellman
protocol, construct multi party procedure for the protocol \cite{AAG}, and analyze the
security of both protocols realized in matrix groups over rings.

The question on finding probabilistic public-key cryptosystems in which the
decryption function has a homomorphic property goes back to~\cite{RAD} (see also~\cite{FM}).
In such a cryptosystem the spaces of messages and of ciphertexts are algebraic
structures $G$ and $H$ and the decryption function $D:G\to H$ is a homomorphism.
A number of such cryptosystems is known for abelian groups, e.g. the quadratic residue
cryptosystem \cite{GB} and its generalization for highest residues \cite{NS} (see also
an overview in \cite{GP}). In most of them the security is based on the intractability of
theoretical number problems close to the integer factoring. Recently, several
homomorphic cryptosystems were constructed for infinite (but finitely presented) groups,
see \cite{GP,GP1} and references there. In this paper we construct one more homomorphic
cryptosystem with $G$ being a free group the trapdoor of which uses a secret permutation
of the generators of $G$.

The third problem considered in this paper is how to produce instances for
cryptosystems based on computations with matrix groups over rings. In contrast to numerous
theoretical cryptosystems where there is a lot of efficient algorithms to generate integers
with given properties (e.g., the pairs of two distinct large primes of the same bit size
used in the quadratic residue cryptosystem), it is not clear a priory how to find efficiently
matrix groups in which some problems (like membership or conjugacy) arising in cryptography
are computationally difficult. We propose a general scheme for solving this problem and
give a specialization of this scheme for matrix groups over finite commutative rings.

In Section~\ref{f010804a} we study key agreement protocols between two parties (named
usually Alice and Bob).
The security of the Diffie-Hellman protocol relies on the difficulty of the following
{\it transporter problem}: having
an action $G\times V\to V$ of a group $G$ on a set $V$ for given $u,v\in V$ to
find $g\in G$ (provided that it does exist) such that $(g,u)\mapsto v$. In case
of $V$ being a cyclic group of order $n$ and $G$ being a group acting on $V$ by taking a power
one arrives to the discrete logarithm
problem (usually, $n$ is taken to be prime). The security of the key agreement protocol
of~\cite{AAG} (see also Subsection~\ref{f270205a}) relies on the difficulty
of the conjugacy problem with respect to a subgroup of $G$. In Subsection~\ref{f270205a}
we extend the construction of \cite{AAG} to {\it multi-party} key agreement protocol.
Then in Subsection~\ref{f050405d} we design another generalization of the
Diffie-Hellman protocol to actions of groups $G$ which satisfy a certain {\it identity}.
Clearly, any abelian group satisfies the identity $aba^{-1}b^{-1}=1$ and more
generally, any solvable group with a fixed length of its derived series satisfies
an appropriate commutator identity. The security of our protocol again relies
on the difficulty of the transporter problem for a suitable action of~$G$.

In Section~\ref{f010804b} we consider homomorphic public key cryptosystems (see e.g.~\cite{GP})
in which the decrypting function (known to Alice) is a group homomorphism $f:G\to H$
where the groups $H,G$ play the roles of the spaces of plain and ciphertext messages
respectively. Usually, the security of a homomorphic cryptosystem relies on the
difficulty of the problem of the membership to a normal subgroup of~$G$ (here, the kernel of~$f$).
Also in Section~\ref{f010804b} we describe a homomorphic cryptosystem in
which as $G$ a free group is taken. This cryptosystem modifies one from~\cite{GP1}
where as $G$ a subgroup of the modular group $\SL_2(\ZZ)$ was considered.
The security of this cryptosystem relies on the difficulty of a certain word
problem. A private key of Alice is an appropriate  permutation of the generators
of the free group $G$, this differs our cryptosystem from the one produced in~\cite{VJS}.

The crucial role in the classical cryptographic constructions (like RSA, discrete
logarithm or quadratic residue~\cite{GB}) plays the natural action of the group
$\aut(\ZZ_n^*)$ on the group $\ZZ_n^*$. So, varying $n$ one gets a mass
pool of instances for cryptographic primitives. This action is a special case of the natural
action of the group $Aut_R(V)$ (viewed as a matrix group) on the free module $V$ over the ring $R$.
In this paper we propose a construction of a pool of matrix groups instances for cryptographic
primitives (Subsection~\ref{f050405u}). The
security of these instances relies on the difficulty of certain problems on matrix groups (e.g.
the membership to a subgroup or the conjugacy with respect to a subgroup). For the complexity
of such problems few results were established in case of matrix groups over fields \cite{BB,KS};
for matrix groups over arbitrary rings much less is known.

The common way in cryptography of producing a trapdoor and a cryptosystem, is to generate a
private key departing from a pair of primes $p,q$, while their product $n=pq$ plays the role of
a public key. In our scheme (see Subsection~\ref{f050405t}) as a private key we take a rooted tree whose
leaves being furnished with specially chosen (non-abelian, in general) groups. We assume that
Alice has in possession such representations of these groups which allow her to solve
efficiently a problem lying in the background  of a cryptosystem (like membership or
conjugacy). Internal
vertices of the tree are endowed with certain operations on groups which allow one to assign
recursively a group to each vertex of the tree starting with its leaves. At the
end of the recursion a group is assigned to the root, and this group plays the
role of a public key. This scheme is also modified to produce a homomorphism of matrix
groups as a public key. In Subsection~\ref{f050405u} we give a realization of this general
scheme in finite matrix groups.

The similarity of the common constructions in cryptography based on commutative groups (say.
$\ZZ_n^*$) with our construction (relying on finite matrix groups) allows us to call the latter
type of constructions the {\it non-commutative cryptography}.

\section{Group theoretical key agreement protocol}
\label{f010804a}

\sbsnt{A multi-party protocol.}\label{f270205a}
 The following
group theoretical variant of key agreement two-party protocol was proposed
in~\cite{AAG}. Let $G$ be a group, and to two parties $A$ and $B$ are assigned
their subgroups
\qtnl{f010804c}
G_A=\lg a_1,\ldots,a_m\rg,\quad G_B=\lg b_1,\ldots,b_n\rg.
\eqtn
The group $G$ and the elements $a_i$, $b_j$ are publicly known.
The parties $A$ and $B$ choose secret elements $a\in G_A$ and $b\in G_B$
and transmit to each other the collections
$$
X_B=\{a^{-1}b_ja\}_{j=1}^n,\quad X_A=\{b^{-1}a_ib\}_{i=1}^m
$$
respectively. Since $A$ (resp. $B$) has a representation of the element $a$ (resp. $b$)
via generators $a_1,\ldots,a_m$ (resp. $b_1,\ldots,b_n$), then $A$ (resp. $B$) can
compute a representation of the element $b^{-1}ab$ (resp. $a^{-1}ba$) via elements of
the set $X_A$ (resp. $X_B$). Thus $A$ and $B$ have a common key
$$
a^{-1}(b^{-1}ab)=[a,b]=(a^{-1}ba)^{-1}b.
$$
An obvious necessary condition for this protocol to be secure is that the set of all such
commutators with $a\in G_A$, and $b\in G_B$ would contain at least two elements.

Let us describe a generalization of the group theoretical key agreement protocol for
$s$ parties with $s\ge 2$ and a single public communicating channel.
Without loss of generality we assume that $s=2^t$ for some
$t\ge 1$, for otherwise in the recursive construction below we divide the parties into
two unequal subsets which leads just to slight changing the notation. As in the case
$s=2$ the groups $G_1,\dots,G_s\subset G$ of the parties are given publically by
their sets of generators. At the initial step the $i$th party chooses a secret key
$a_i\in G_i$, $i=1,\ldots,s$. Let $S_1$ and $S_2$ be disjoint $s/2$-subsets of
the set $\{1,\ldots,s\}$. Then given $u=1,2$ the parties from $S_u$ recursively construct
the common key $K_u\in G$, such that for all $i\in S_u$ there exist integers
$\varepsilon_{i,j}\in\{-1,+1\}$ and $1\le m_i\le s/2$, and certain elements
$B_{i,1},\ldots,B_{i,m_i}\in\lg \{a_j:\ j\in S_{u,i}\}\rg$ with
$S_{u,i}=S_u\setminus\{i\}$, for which we have
$$
K_u=(B_{i,1}^{-1}a_i^{\varepsilon_{i,1_{}}}B^{}_{i,1})\cdots
(B_{i,m_i}^{-1}a_i^{\varepsilon_{i,m_i}}B^{}_{i,m_i}).
$$
By recursion we can assume that the $i$th party knows the elements
$B_{i,j}^{-1}aB_{i,j}^{}$ for all $j$ and for all chosen generators $a$ of the
group $G_i$ (and thereby, it knows $B_{i,j}^{-1}a_iB_{i,j}^{}$), but does not necessary
know $B_{i,j}$. At this point the party $i\in S_u$ sends the elements
$B_{i,j}^{-1}aB_{i,j}^{}$ for all the chosen generators $a$ of the group $G_i$ to a
certain party from the set $S_{u'}$ with $u'=3-u$ and asks for the elements
$K_{u'}^{-1}B_{i,j}^{-1}aB_{i,j}^{}K_{u'}^{}$. Then for $u=1$ the $i$th party
computes the element
$$
[K_1,K_2]=K_1^{-1}(K_2^{-1}K_1^{}K^{}_2)=
K_1^{-1}(K_2^{-1}(B_{i,_{m_i}}^{-1}a_i^{\varepsilon_{i,m_i}}B^{}_{i,m_i})K_2^{})\cdots
(K_2^{-1}(B_{i,1}^{-1}a_i^{\varepsilon_{i,1}}B^{}_{i,1})K_2^{}).
$$
Similarly, for $u=2$ the $i$th party computes the element
$[K_1,K_2]=(K_1^{-1}K_2^{}K_1^{})^{-1}K_2$. Thus this element
can be chosen as the common key for all parties. It is easy to see that the $i$th
party computes the common key in $O(s|a_i|)$ operations in the group $G$, where $|a_i|$
denotes the length of the word $a_i$ in the chosen generators of the group $G_i$.

\sbsnt{A new protocol.}\label{f050405d}
In this subsection we define a new group-theoretical two party key agreement
protocol that can be viewed as a non-commutative generalization of the Diffie-Hellman
protocol (see~\cite{GB}).

Let $G$ be a group acting on a set $X$ so that given $(x,g)\in X\times G$ the
image $x^g$ of $x$ with respect to $g$ can be efficiently computed. Two parties
$A$ and $B$ going to choose a secret common key from $X$, fix publically
subgroups $G_A,G_B$ of the group $G$ and two words
$$
W_A(u_A,u_B)=u_A^{a_{1,1}}u_B^{b_{1,1}}\cdots u_A^{a_{1,m_1}},\qquad
W_B(u_A,u_B)=u_B^{b_{2,1}}u_A^{a_{2,1}}\cdots u_B^{b_{2,m_2}}
$$
of the free group $F_2$ with two free generators $u_A,u_B$ such that
\nmrt
\tm{W1} $m_1,m_2\in\N$, $a_{i,j},b_{i,j}\in\ZZ$ for all $i,j$, and $a_{1,m_1}\ne 0$,
$b_{2,m_2}\ne 0$,
\tm{W2} $W_A(g_A,g_B)=W_B(g_A,g_B)$ for all $(g_A,g_B)\in G_A\times G_B$.
\enmrt
The protocol begins with the choice of a publically known element $x_0\in X$ and
the secret elements $g_A\in G_A$ by the party $A$ and $g_B\in G_B$ by the party $B$.
Then during the communications the party $A$ performs the following:
\nmrt
\item[--] Set $K_A=x_0$.
\item[--] For $i=1,\ldots,m_1-1$ send $K_A^{{g_A}^{a_{1,i}}}$ and receive $K_A:=K_A^{{g_A}^{a_{1,i}}{g_B}^{b_{1,i}}}$.
\item[--] Set $K_A:=K_A^{{g_A}^{a_{1,m_1}}}$.
\enmrt
The communications of the party $B$ are defined similarly. Thus at the end of the
communication process due to condition (W2) the parties $A$ and $B$ have the common key
$$
K_A=x_0^{W_A(g_A,g_B)}=x_0^{W_B(g_A,g_B)}=K_B.
$$
For $X=\ZZ_p^*$ with $p$ being a prime, $G=G_A=G_B$ being the group
$\ZZ_{p-1}^*\cong\aut(\ZZ_p^*)$
and $W_A(u_A,u_B)=u_Bu_A$, $W_B(u_A,u_B)=u_Au_B$ we come to the Diffie-Hellman
protocol.

This scheme can be easily realized for a solvable group $G$ with bounded length $n$
of the derived series of $G$. For example, one can take $G_A=G_B=G$ and choose the words
$W_A=W_{A,n}$ and
$W_B=W_{B,n}$ by induction on $n$ as follows. If $n=1$, then the group $G$ is abelian
and so conditions (W1) and (W2) are satisfied for
$$
W_{A,1}(u_A,u_B)=u_Bu_A,\qquad
W_{B,1}(u_A,u_B)=u_Au_B.
$$
For $n\ge 2$ the commutator $[g,h]=g^{-1}h^{-1}gh$ with arbitrary $g,h\in G$ belongs to
the derived subgroup $G'=[G,G]$ of $G$ (the derived length of $G'$ equals $n-1$).
Assume by induction that conditions (W1) and (W2) are satisfied for the words
$W_{A,n-1}$ and $W_{B,n-1}$. Then a straightforward checking shows that these conditions
are also satisfied, for example, for the words
$$
W_{A,n}=W_{A,n-1}([u_B^{},u_A^{}],[u_A^{-1},u_B^{-1}]),\qquad
W_{B,n}=W_{B,n-1}([u_B^{},u_A^{}],[u_A^{-1},u_B^{-1}]).
$$
This follows from the fact that the length (the number of letters) of the word $W_{A,n}$
(as well as $W_{B,n}$) equals $2\cdot 4^{n-1}$ which one can verify by induction on~$n\ge 1$.
More generally, to define $W_{A,n}$ and $W_{B,n}$ one can choose arbitrary words
$W_1,W_2,W_3,W_4\in W_X$ where $X=\{u_A,u_B\}$ and $W_X$ is the set of all words in the
alphabet $X^\pm$, and use $[W_1,W_2]$ and $[W_3,W_4]$ instead of $[u_A^{},u_B^{}]$
and $[u_B^{-1},u_A^{-1}]$ respectively. Certainly, to provide condition (1) one should
guarantee that the words $W_{A,n-1}(u_A,u_B)$ (resp. $W_{B,n-1}(u_A,u_B)$) and $W_2$
(res. $W_4$) must be terminated to $u_A$ (resp. $u_B$). To avoid triviality we also
should take $W_1,\ldots,W_4$ so that $W_{A,n}$ and $W_{B,n}$ would be nonidentity
elements in the underlying free group.

Clearly, any realization of the above protocol is based on identities of the
group $G$. In addition to commutator identities for solvable groups (see above) one can also
use the identity $x^m=1$ (that holds in any finite group the order of which is
a divisor of $m$, and in the Burnside groups). In this case we can choose as the words
$W_A$ and $W_B$ the prefix and the inverse of the suffix of the word $(u_Au_B)^m$,
respectively, so that
the prefix is terminated to $u_A$. In fact, as it was
proved by B.Neumann any variety of groups can be given by a collection of identities such
that the first of them is of the form $x^m=1$ with $m$ being a nonnegative integer,
whereas the other ones are the elements of the commutant of the underlying free
group (see~\cite{MKS}).

We complete the subsection by making two remarks on the above protocol. First,
the set $X$ must be of superpolynomial size, for otherwise the key agreement scheme can be broken
in polynomial time by the
known permutation group theory technique (see~\cite{L}). Second, the words $W_A$
and $W_B$ must be chosen so that the number of elements $W_A(g_A,g_B)=W_B(g_A,g_B)$
with $g_A,g_B\in G$ would contain at least two elements.

\sbsnt{On the security of the protocols.}\label{f050404c}
In the above protocols we assume that all groups are given explicitly, e.g. by
sets of generators, so that the group operations can be performed efficiently. Then
the security of the first protocol is based on the intractability of the following
problem (see~\cite{SU}).

\prblm{Subgroup Conjugation Search Problem (SCSP)}{Given a group $G$, subgroups $H_1,H_2$
of $G$, and two elements $f,g\in H_1$, find an element $h\in H_2$ such that
$f = h^{-1}gh$, provided that at least one such $h$ exists.}

As usually in the cryptography, an efficient algorithm solving SCSP would break the
protocol (but to break the protocol it is not necessary to solve SCSP). Such an algorithm
does exist for $G=\GL(n,\F_q)$ where $n$ is a natural number, $\F_q$ is a finite field of
the order~$q$, and the subalgebra $A(H_2)$ of the full matrix algebra $\Mat_n(\F_q)$ generated
by the group $H_2$ is such that
$$
A(H_2)\cap G=H_2.
$$
Then for arbitrary $H_1$ the problem SCSP can be solved in probabilistic polynomial time
(in $n$ and in $\log q$) by the linear algebra technique, provided that $n$ is less than $q/2$.
Indeed, in this case the solution of the linear system $hf-gh=0$ with respect to $h\in A(H_2)$
is an element of $H_2$ with a great probability. (From \cite{CY} it follows that in this case
the problem SCSP can be solved efficiently even by a deterministic algorithm.)

It seems that the problem SCSP remains difficult when $G$ is restricted to
subgroups of the group $\GL(V,R)$ of all invertible $R$-linear transformations
of the free $R$-module~$V$ where $R$ is a finite commutative ring. To see this
we consider the Linear Transporter Problem on the intractability of which the second
protocol is based.

\prblm{Linear Transporter Problem (LTP)}{Let $R$ be a commutative ring, $V$ be an
$R$-module and $G\le\GL(V,R)$. Given $u\in V$ and $v\in u^G=\{u^g:\ g\in G\}$
find $g\in G$ such that $v=u^g$.}

A special case of (LTP) is the Discrete Logarithm Problem. Indeed, take $V=\ZZ_p^*$
with $p$ being a prime. Then $V$ can be considered as an one-dimensional module over the ring
$R=\End(V)\cong\ZZ_{p-1}$ (with respect to taking the power
$v\mapsto v^n$ where $v\in V$, $n\in\ZZ_{p-1}$). Choosing $u$ to be a generator of the group
$V$ we come to the Discrete Logarithm Problem.

Preserving the notation of LTP set $T(V)=\{T_v:x\mapsto x+v,\ v,x\in V\}$ to be the
translation group of the $R$-module $V$. Then obviously
$$
v=u^g\ \Leftrightarrow\ T_v=g^{-1}T_ug,\qquad u,v\in V,\quad g\in\GL(V,R).
$$
So the problem LTP is the special case of the problem SCSP with $G=\AGL(V,R)$, $H_1=T(V)$
and $H_2=\GL(V,R)$. (Here $\AGL(V,R)=T(V)\GL(V,R)$ is the group of all affine transformations
of $V$.) This shows that
SCSP is at least as hard as LTP. In particular, this construction gives us a family
of groups for which the problem SCSP turns to be at least as hard as the Discrete
Logarithm Problem. A general technique to construct groups of this kind will be
given in Section~\ref{f3107a}.

\section{Homomorphic cryptosystems over groups}\label{f010804b}

\sbsnt{A general scheme.}
A homomorphic cryptosystem is a probabilistic public key scheme (in the sense of
\cite{GB}) in which the spaces of plaintext messages and ciphertexts are groups $H_k$ and
$G_k$ respectively, depending on a security parameter $k$ and such that its decryption function
\qtnl{f010804e}
f_k:G_k\to H_k
\eqtn
is an epimorphism for all~$k$. Usually, in a homomorphic cryptosystem the public key
includes generator sets $X_k$ and $Y_k$ of the groups $G_k$ and $H_k$, and some
set $R_k\subset X_k$ such that $Y_k\subset f_k(R_k)=\{f_k(g):\ g\in R_k\}$. Besides,
it is assumed that there are publically known $k^{O(1)}$-algorithms to solve the
following problems:
\nmrt
\tm{1} given two elements $a,b$ of $G_k$ (resp. $H_k$) find the element $ab^{-1}$,
\tm{2} given $y\in Y_k$ find an element of the set $R_k\cap f_k^{-1}(y)$,
\tm{3} generate a random element of the group $\ker(f_k)$
\enmrt
where sizes of all elements are assumed to be at most $k$.
Under these assumptions the encryption can be performed in time $k^{O(1)}$ as follows.
First, given a message $h=y_1\cdots y_m\in H_k$ with $y_i\in Y_k$ and $m$ being a
natural number at most $k^{O(1)}$, Bob computes in time polynomial in $k$ an element
$r=r_1\cdots r_m\in G_k$ such that $r_i\in R_k$ and $f_k(r_i)=y_i$ for all $i$.
Second, Bob mixes $r$ with random elements $g_1,\cdots,g_{m+1}\in G_k$ belonging to
the kernel of the homomorphism $f_k$ and outputs the element
$g=g_1r_1g_2\cdots g_mr_mg_{m+1}$ as the ciphertext of $h$. Alice
being able to compute $f_k$ efficiently performs the decoding as follows:
$$
f_k(g)=f_k(g_1r_1g_2\cdots g_mr_mg_{m+1})=
f_k(r_1)\cdots f_k(r_m)=y_1\cdots y_m=h.
$$
The key point of such a system is to choose a presentation of the group~$G_k$ and
the epimorphism $f_k$ in order to provide the inverse of $f_k$ to be a trapdoor
function. The exact definition of homomorphic public-key cryptosystems and a survey of
constructions can be found in~\cite{GP,GP1}.

One way to implement the general concept of a homomorphic cryptosystem is to
take $G_k$ to be a subgroup of a certain group $F$ such that the group operations in $F$
can be performed in time polynomial in the size of operands. In the cryptosystems
from \cite{GP} and \cite{GP1} the group $F$ was taken as a free product of abelian
groups and a modular group, respectively. In these cryptosystems the restriction of
the mapping $f_k$ to the set $R_k$ was known publically and one can produce efficiently
random $k^{O(1)}$-size elements of the group $\ker(f_k)$. In fact, the security of
these cryptosystems was based on the difficulty of the membership problem (see below)
for special subgroups of the group $G_k$. In the next subsection we present a new
homomorphic public-key cryptosystem of this kind (but with another trapdoor).

\sbsnt{A new homomorphic scheme.}
Let $H=\lg Y;\R\rg$ be a finitely presented group generated by the set $Y$ of
cardinality $k\ge 2$ with $\R\subset W_Y$ as the set of relations. As the group $F$
mentioned above we take the free group $\lg Y\rg$. For a permutation $\sigma\in\sym(Y)$
denote by $\varphi_\sigma$ the automorphism of the group $F$ induced by~$\sigma$. Set
\qtnl{f270305}
X=X_\sigma=\{\varphi_\sigma^{-1}(r_yyr'_y):\ y\in Y\}
\eqtn
where $r_y$ and $r'_y$ are randomly chosen words of size $O(k)$ belonging to the set
$W_\R\subset W_Y$. Then $G=\lg X\rg$ is a subgroup of the group $F$. Moreover,
the mapping $f_\sigma:G\to H$ defined by a commutative diagram
\qtnl{f240305}
\begin{array}{ccc}
G                                     & \stackrel{f_\sigma}{\longrightarrow} &     H\\
\Big\downarrow\vcenter{\llap{$\id_G$\phantom{\ \ }}}  & & \Big\uparrow\vcenter{\rlap{$\rho$}}\\
F                                     & \stackrel{\varphi_\sigma}{\longrightarrow} & F\\
\end{array}
\eqtn
where $\rho:F\to H$ is the epimorphism induced by the mapping $\id_Y$, is an epimorphism such
that given $x\in X$ we have (see (\ref{f270305})):
$$
f_\sigma(x)=\rho(\varphi_\sigma(x))=\rho(\varphi_\sigma(\varphi_\sigma^{-1}(r_yyr'_y)))=
\rho(r_yyr'_y)=\rho(r_y)\rho(y)\rho(r'_y)=\rho(y)=y
$$
where $y$ is the element of $Y$ for which $x=\varphi_\sigma^{-1}(r_yyr'_y)$ (see~(\ref{f270305})).
In particular, $f_\sigma(X)=Y$ and the restriction of $f_\sigma$ to $X$ is a
bijection. This enables us to construct a homomorphic cryptosystem as follows.

\vspace{4mm}
\noindent{\bf Secret Key:} the permutation $\sigma\in\sym(Y)$.
\vspace{2mm}

\noindent{\bf Public Key:} a natural number $k\ge 2$, a group $H=\lg Y;\R\rg$ with $|Y|=k$,
a subgroup $G=\lg X_\sigma\rg$ of the free group $F=\lg Y\rg$, and the bijection
$f:X_\sigma\to Y$ coinciding with the restriction of the homomorphism $f_\sigma$ to $X_\sigma$.
\vspace{2mm}

\noindent{\bf Encryption:} a message $M=y_{i_1}\cdots y_{i_t}\in H$ where $y_{i_j}\in Y^\pm$,
is encrypted by the element
$$
E(M)=f^{-1}(s_1y_{i_1}s'_1)\cdots f^{-1}(s_ty_{i_t}s'_t)\in G
$$
where $s_i$ and $s'_i$ are random words of the set $W_\R\subset W_Y$ of size $O(k)$, and for a word
$w=\cdots y \cdots\in W_Y$ we set $f^{-1}(w)=\cdots f^{-1}(y)\cdots$.
\vspace{2mm}

\noindent{\bf Decryption:} a ciphertext $C=y_{i_1}\cdots y_{i_t}\in G\subset F$ where
$y_{i_j}\in Y^\pm$, is decrypted to $D(C)=y_{i_1}^\sigma\cdots y_{i_t}^\sigma\in H$.
\vspace{4mm}

To prove the correctness of the decryption we note that
$f=(f_\sigma)|_{_{X_\sigma}}$, $f_\sigma=\rho(\varphi_\sigma)|_G$, and
$\varphi_\sigma(y)=y^\sigma$ for all $y\in Y$ (see~(\ref{f240305})). Since obviously
$\varphi_{\sigma}^{-1}=\varphi_{\sigma^{-1}}$, we have
$$
D(E(y_{i_1}\cdots y_{i_t}))=
D(f^{-1}(s_1y_{i_1}s'_1)\cdots f^{-1}(s_ty_{i_t}s'_t))=
D(\varphi_{\sigma^{-1}}(s_1y_{i_1}s'_1\cdots s_ty_{i_t}s'_t))=
$$
$$
\varphi_{\sigma^{-1}}(s_1y_{i_1}s'_1)^\sigma\cdots\varphi_{\sigma^{-1}}(s_ty_{i_t}s'_t)^\sigma=
s_1y_{i_1}s'_1\cdots s_ty_{i_t}s'_t=y_{i_1}\cdots y_{i_t}.
$$
Clearly, that both encryption and decryption algorithms are polynomial-time in the
size of the input words.

The security of the homomorphic cryptosystem will be discussed in the
next subsection. Here we only make several remarks on the possible implementations.
First, we note that
it is not necessary to work with words; instead
of this one can use a matrix representation of the group $F$ (see~\cite{GP1}). Next, to choose
the set $Y$ so that $|Y|=k$, one can take any set $S$ of generators of $H$ and add
to it $k-|S|$ random  elements of $H$ whenever $|S|<k$. Finally, as in
Section~\ref{f010804a} any implementation of the above cryptosystem must be supported
by sufficiently large class of candidates for groups $H$. We will return to this
problem in Section~\ref{f3107a}.

\sbsnt{On the security of homomorphic schemes.}\label{f090405c}
Concerning the security of the homomorphic cryptosystem suppose first that the order of
the group $H$ is at most $k^{O(1)}$ (e.g. such an assumption was done in~\cite{GP}).
Then using the generator set $Y$ of $H$ one can list all the elements $h_1,\ldots,h_m$
of this group in time $k^{O(1)}$ and then to find within the same time a set $\{g_1,\ldots,g_m\}$
of distinct representatives of right cosets of $G$ by $G_\sigma=\ker(f_\sigma)$ (one can set
$g_i=f^{-1}(h_i)$ for all $i$). Now if an adversary Charlie could recognize efficiently the elements of $G$
belonging to $G_\sigma$, then he would efficiently compute $f_\sigma(g)$ for all $g\in G$
due to the formulae
$$
f_\sigma(g)=f_\sigma(g_i)\ \Leftrightarrow\ gg_i^{-1}\in G_\sigma
$$
where $i\in\{1,\ldots,m\}$. Thus in this case the security of our cryptosystem is based on the
intractability of the following problem:
\vspace{2mm}

\prblm{Membership Testing (MT)} {Given a group $F$ and its subgroup $G$ test whether a given
$g\in F$ belongs to~$G$.}

Suppose now that the order of $H$ to be arbitrary. Then a quite natural way to break the cryptosystem
is to find an expression of any $g\in G$ in the terms of generators belonging to the set
$X_\sigma$ (the attack of this kind was considered in~\cite{GP1}). Indeed, if Charlie could
find efficiently for any element $g\in G$ an expression $g=x_1\cdots x_m$ where
$x_i\in X_\sigma^\pm$ for all $i$, then he would efficiently compute $f_\sigma(g)$ due to formulae
$$
f_\sigma(g)=f_\sigma(x_1)\cdots f_\sigma(x_m)=f(x_1)\cdots f(x_m)
$$
(we recall that the bijection $f:X_\sigma\to Y$ is given publically). Thus in this case we
come to the presentation problem (see~\cite{GP1}). The MT problem and the presentation problem
are closely related each to other (but generally could be not polynomial-time equivalent) and
one can combine them in the following well-known problem of computational group theory
(see \cite{BB}).
\vspace{2mm}

\prblm{Constructive Membership Testing (CMT)} { Given a group $F$ and its subgroup $G$
generated by a set $X$ find an expression of a given $g\in F$ as a word in $X$,
or determine that $g\not\in G$.}

Last two decades a great attention was paid to CMT with different presentations of
the group $G$. For example, if $F$ is a subgroup of the symmetric group of degree~$n\ge 1$,
then the CMT can be solved in time $n^{O(1)}$ by the sift algorithm (see e.g.~\cite{L}).
In the case of groups $F=\GL(n,\F)$ where $\F$ is an algebraic number field,
there exists an effective Las Vegas algorithm solving CMT \cite{BB}. However, for $n=1$
and $\F$ being a finite field, CMT is nothing else but the the Discrete Logarithm Problem.
In \cite{BB} it was conjectured that CMT is difficult whenever the group $G$ either
involves a large abelian group as a quotient of a normal subgroup or has nonabelian
composition factors which require large degree permutation representations. Finally,
the problem becomes much more difficult if we take $F=\GL(n,R)$ the group of $n\times n$
invertible matrices over a ring $R$. In this case the problem is undecidable for $n=4$
and $R=\ZZ$ (see~\cite{M58}).

\section{Cryptographical generation of groups}\label{f3107a}

\sbsnt{A general scheme.}\label{f050405t}
We begin with a general scheme to construct a vast family of groups and homomorphisms
supporting both key agreement protocols of Section~\ref{f010804a} and
homomorphic cryptosystems of Section~\ref{f010804b}. Let $\G$ be a class of groups
closed with respect to a set $\O$ of group-theoretical operations of different
arities (like direct or wreath products). For an integer $s\ge 1$ we denote by $\O_s$
a set of all operations of arity $s$ belonging to~$\O$. For a set $\G_0\subset\G$ we
define recursively a class $\P(\G_0,\O)$ of pairs $(G,T)$ where $G\in\G$ and $T$ is a
rooted labeled tree, as follows:
\vspace{2mm}

\noindent{\bf Base of recursion:} any pair $(G,T)$ with $G\in\G_0$ and $T$ being the one-point tree with root
labeled by $G$, belongs to $\P(\G_0,\O)$.
\vspace{2mm}

\noindent{\bf Recursive step:} given pairs $(G_1,T_1),\ldots,(G_s,T_s)\in\P(\G_0,\O)$ and an operation $o\in\O_s$,
the class $\P(\G_0,\O)$ contains the pair $(G,T)$ where
$G=o(G_1,\ldots,G_s)$ and $T$ is the tree obtained from $T_1,\ldots,T_s$ by
adding a new root labeled by $o$ and the sons being the roots of $T_1,\ldots,T_s$.
\vspace{2mm}

Let $(G,T)\in\P(\G_0,\O)$. Then obviously $G\in\G$ and the {\it derivation tree} $T$ of
$G$ provides the constructive proof for this membership. The group $G$ is uniquely
determined by $T$ and we call it the {\it group associated with $T$}. The fact, that a derivation
tree is an ordinary rooted tree the leaves and the internal vertices of which are
labeled by elements of $\G_0$ and $\O$ respectively, enables us to choose a random
derivation tree of a fixed size.

Suppose from now on that all the groups of $\G$ are given in a certain way (e.g., one
can take as $\G$ a class of matrix groups given by generator sets). We assume also
that for each operation $o\in\O_s$ and groups
$G_1,\ldots,G_s\in\G$, the size $L(G)$ of the
presentation of the group $G=o(G_1,\ldots,G_s)$ is at most $O(L)$ where
$L=\sum_{i=1}^sL(G_i)$ and the group $G$ can be constructed from $G_1,\ldots,G_s$ in
time $L^{O(1)}$. Let us define a size $L(T)$ of a derivation tree $T$ to be the sum of the sizes
of all labels of $T$; thus $L(T)$ includes the sizes of the groups assigned to the leaves
of $T$ together with the number of edges of $T$. Then the size of any pair
$(G,T)\in\P(\G_0,\O)$ is $O(L(T))$, and the knowledge of $T$ enables us to find $G$ in
time polynomial in $L(T)$.

One of the problems arising in constructions of group-theoretical public key cryptosystems
is to find an efficient algorithm to produce a random group (or a collection of groups)
belonging to a special class $\G$ and with a given size $L$ of the presentation.
Such a group $G$ must be equipped with a private key providing an efficient solution of
a certain computational problem for $G$ that is supposedly difficult in the class $\G$
without knowledge of a private key.
Our approach to the above problem is to choose an appropriate class $\G_0$ of groups,
a set $\O$ of group-theoretical operations, and then to generate instances for the
cryptosystem in question as follows:
\vspace{2mm}

\noindent{\bf Step 1:}
given a security parameter $L$ choose randomly groups $G_1,\ldots,G_t\in\G_0$,
such that $\sum_{i=1}^tL(G_i)=O(L)$;
\vspace{2mm}

\noindent{\bf Step 2:}
choose randomly a rooted labeled tree $T$ of size $O(L)$ and with $t$ leaves
being labeled by $G_1,\ldots,G_t$;
\vspace{2mm}

\noindent{\bf Step 3:}
compute the group $G$ associated with $T$ (i.e. $(G,T)\in\P(\G_0,\O)$);
\vspace{2mm}

\noindent{\bf Step 4:}
output the group $G$ as a public key and the labeled tree $T$ as a secret key.
\vspace{2mm}

Denote by $\G^*$ the class of groups $G$ such that $(G,T)\in\P(\G_0,\O)$ for some labeled
tree $T$. Then the secrecy of the key $T$ is based on the intractability of the
following problem: given $G\in\G^*$
find a derivation tree $T$ associated with $G$. A special case of this problem will be
considered in Section~\ref{f090405a}.

For a homomorphic cryptosystem the above scheme is not sufficient because together
with the group $G$ we have to provide a group $H$ and a secret homomorphism $f:G\to H$.
To this end suppose that each group $G\in\G_0$ is equipped with a set $M(G)$ of pairs
$(H,f)$ where $H\in\G_0$ and $f:G\to H$ is a homomorphism. We also assume that
given homomorphisms $f_i:G_i\to H_i$ with $G_i,H_i\in\G^*$ for $i=1,\ldots,s$, and an
operation $o\in\O_s$ there exists an efficiently computed homomorphism $f:G\to H$ where
$G=o(G_1,\ldots,G_s)$ and $H=o(H_1,\ldots,H_s)$ such that $f|_{G_i}=f_i$ for all~$i$
(here we suppose in addition that $G_i$ is a subgroup of $G$). This homomorphism is denoted
by $o(f_1,\ldots,f_s)$. In this notation the set $\M(\G_0,\O)$ of
instances $(G,f)$ for a homomorphic cryptosystem can be defined recursively as follows:
\vspace{2mm}

\noindent{\bf Base of recursion:} any pair $(G,f)$ with $G\in\G_0$ and $f\in M(G)$ belongs
to the set $\M(\G_0,\O)$;
\vspace{2mm}

\noindent{\bf Recursion step:} given pairs $(G_1,f_1),\ldots,(G_s,f_s)\in\M(\G_0,\O)$ and an
operation $o\in\O_s$, the class $\M(\G_0,\O)$ contains the pair $(G,f)$ where
$G=o(G_1,\ldots,G_s)$ and $f=o(f_1,\ldots,f_s)$.
\vspace{2mm}

We observe, that in the process of constructing the homomorphism $f:G\to H$ we also
produce the derivation trees of the groups $G$ and $H$. A realization
of these general schemes in finite matrix groups will be considered in the next subsection.

\sbsnt{Generating matrix groups.}\label{f050405u}
Let us define the classes $\G_0,\G$ of groups and the set $\O$ of operations. First, we set
$$
\G=\cup_n\cup_R \{G:\ G\ \mbox{is a subgroup of}\ \GL(n,R)\}
$$
where $n$ and $R$ run over natural numbers and finite commutative rings respectively. Thus
any $G\in\G$ is a group of $n\times n$ invertible matrices with entries belonging to $R$
for some $n\in\N$ and some finite commutative ring $R$. We recall that any such ring is
a direct sum of local commutative rings and each of the latter can be described via
appropriate Galois ring: the Galois ring $\GR(p^m,r)$ of characteristic $p^m$ and rank $r$ is
$\ZZ_{p^m}[x]/(f)$ where $f\in\ZZ_{p^m}[x]$ is a monic polynomial of degree $r$ whose image in
$\ZZ_p[x]$ is irreducible (see \cite{MD74}).

\prpstnl{f090505a}{\rm \cite{MD74,Wan}}
Let $R$ be a finite commutative local ring of characteristic $p^m$ and $\F=\GF(p^r)$ the residue
field of $R$. Then
\nmrt
\tm{1} $R^\times=\T\times(1_R+\rad(R))$ where $\T$ is a cyclic group isomorphic to $\F^\times$,
\tm{2} the subring $R_0$ of $R$ generated by $\T$ is a Galois ring $\GR(p^m,r)$,
\tm{3} $R$ is a homomorphic image of the ring $R_0[X_1,\ldots,X_t]$ where $t$ is the minimal
size of a generator set of the radical of $R$.\bull
\enmrt
\eprpstn

\prpstnl{f090505b}{\rm \cite{MD74}}
Let $p$ be a prime and $m,r$ be natural numbers. Then
\nmrt
\tm{1} there exists the unique up to isomorphism Galois ring $\GR(p^m,r)$ of characteristic
$p^m$ and rank $r$,
\tm{2} each element $r\in\GR(p^m,r)$ is uniquely represented in the form $r=\sum_{i=0}^{m-1}t_ip^i$
where $t_i\in\T\cup\{0\}$ for all~$i$,
\tm{3} given $\ov\sigma\in\aut(\F)$ the mapping $r\mapsto\sum_{i=0}^{m-1}t_i^\sigma p^i$
where $\sigma$ is the automorphism of the group $\T$ induced by $\ov\sigma$ (see
statement (1) of Proposition~\ref{f090505a}), is an automorphism of~$\GR(p^m,r)$.\bull
\enmrt
\eprpstn

Due to statements (2),(3) of Proposition~\ref{f090505a} and statements (2) of
Proposition~\ref{f090505b} a representation of the finite
commutative ring $R$ (resp., the group $G$) can be chosen to be polynomial in $\log(|R|)$ (resp.
in $n$ and $\log(|R|)$). We
also admit a hidden representation of $R$ in which the decomposition in local summands is not
presented explicitly, for example the ring of residues modulo an integer can be completely given
by indicating this integer.

We define a set $\G_0\subset\G$ to be a class of classical simple (including abelian) subgroups
$G$ of the groups $\GL(n,\F)$ where $n\in\N$ and $\F$ is a finite field. Any such group
$G\in\G_0$ is given by a set of generators so that the Membership Testing Problem can be solved
in time polynomial in $n$ and in the bit size of $\F$. (Indeed, any nonabelian classical matrix
group can be given together with a suitable matrix representation which can be used for testing
membership; for an abelian group of a prime order~$p$ one can use, e.g. the two-dimensional
representation
\qtnl{f210405a}
\ZZ_p^+\to\GL(2,p),\quad x\mapsto
\begin{pmatrix}
 1 & x \\
 0 & 1 \\
\end{pmatrix}
\eqtn
which gives a trivial membership testing algorithm).
In fact, it is not necessary that $\G_0$ contains all classical groups; one can form
$\G_0$ from the group of special types, e.g. $\PSL(n,\F)$ or something like that. Since
the elements of $\G_0$ are parametrized by the tuples of naturals, one can efficiently choose
a random group $G\in\G_0$ with a given size $L(G)$ of presentation.

The choice of the set $\O$ of operations was inspired by the Aschbacher theorem \cite{Asch}
on classifying maximal subgroups of classical groups. Let us describe the operations.
\vspace{2mm}

{\bf Changing the underlying ring.} Let $R$ be a finite commutative ring and $R'$ be an
extension of $R$. Then the natural monomorphism
$$
\varphi:\GL(n,R)\to\GL(n,R')
$$
gives an unary operation in $\G$ taking $G\in\G$ to $\varphi(G)$. This operation can be
performed efficiently whenever e.g. the embedding $R$ to $R'$ is given explicitly and the
number $d=[(R')^+:R^+]$ is small. Another example is the extension of $\ZZ_m$ to
$\ZZ_{m'}$ where $m$ is a divisor of $m'$.
Conversely, any embedding of the ring $R'$ into the ring $\Mat(d,R)$ induces the natural
monomorphism
$$
\varphi':\GL(n,R')\to\GL(nd,R)
$$
taking a matrix of $\GL(n,R')$ to the block matrix of $\GL(nd,R)$ with $d^2$ blocks
of size $n$. (Such a situation arises e.g. when $R'$ is a field of the order $q^d$ and
$R$ is its subfield of the order $q$, or when $R'$ is isomorphic to the direct sum of $d$ copies
of $R$.) This produces another unary operation in $\G$ taking
$G\in\G$ to $\varphi'(G)$. In order not to blow up the representation one should assume that
$d$ is small. In both cases the isomorphism type of the group $G$ (as an abstract group) does not
change, but the operations change it as a linear group. In fact, our constructions start with
matrix groups over a finite field $\F$. To pass to rings one can use standard extensions
with $R=\F$ and $R'=\Mat(m,R)$, and also with $R=\Mat(n,p)$ and $R'=\Mat(m,\ZZ_{p^d})$ with
a prime $p$.
\vspace{2mm}

{\bf Direct products.}
Suppose that groups $G_1,\ldots,G_s\in\G$ are such that $G_i\le\GL(n_i,R)$ where $n_i\in\N$ and
$R$ is a finite commutative ring. Then
$$
G=G_1\otimes\cdots\otimes G_s\le\GL(n,R)
$$
where $n=\prod_{i=1}^sn_i$, and we obtain an $s$-ary operation in $\G$. A set of generators for
the group $G$ can be efficiently constructed from the generating sets for $G_1,\ldots,G_s$
by means of the Kronecker product of the corresponding matrices. When $R$ is a field the group
$G$ is irreducible iff so are the groups $G_1,\ldots,G_s$. (A matrix group $G$ is called
irreducible if the underlying linear space contains no nontrivial $G$-invariant subspaces.)

Similarly, if $m=n_i$, $G_i\cap G_i'=\{I_m\}$ and $G'_i$ normalizes $G_i$ where
$i=1,\ldots,s$ and $G'_i$ is the group generated by $G_j$, $j\ne i$, then
$G_1\times\cdots\times G_s$ is a subgroup of $\GL(m,R)$ which gives one more $m$-ary operation.
\vspace{2mm}

{\bf Wreath products.}
The {\it wreath product} $G\wr\Gamma$ of a group $G$ and a permutation group $\Gamma\le\sym(m)$
is defined to be the semidirect product of the $m$-fold direct product $G^m=G\times\cdots\times G$
by the group $\Gamma$ acting on $G^m$ via coordinatewise permutations. If $G\le\GL(n,R)$, then
the group $G\wr\Gamma$ has two natural linear representations obtained from the natural
monomorphisms
$$
G^m\to\GL(nm,R),\qquad G^m\to\GL(n^m,R),
$$
the first of which is induced by the $m$-fold direct sum of the underlying $R$-module, whereas
the second one is induced by the $m$-fold tensor product of it. The images of the group $G\wr\Gamma$
are called the {\it imprimitive} and the {\it product actions} of the wreath product, respectively.
Thus we obtain two more efficiently computable $m$-ary operations in $\G$. In the case of $R$ being
a field the resulting groups are always irreducible whenever $G$ is irreducible and $\Gamma$ is
transitive.
For our purpose it is enough to set
$\Gamma$ to be the symmetric group. More elaborated way could be based on the
fact that any transitive group is obtained from the action of a group on the set of
right cosets by some subgroup by means of right multiplications.
\vspace{2mm}

{\bf Conjugations.}
An obvious unary operation in $\G$ consists in the conjugation of a group $G\in\GL(n,R)$ by
means of a randomly chosen matrix from $\GL(n,R)$. Such an operation enables us to hide
the form of a generator set of the group $G$.
\vspace{2mm}

Let $\O$ be the set of the above operations and $\G^*\subset\G$ be the set of all groups
$G$ such that $(G,T)\in\P(\G_0,\O)$ for some rooted labeled tree $T$ (see
Subsection~\ref{f050405t}). In the following statement we consider the specializations
of the problems MT (see Subsection~\ref{f090405c}) and LTP (see Subsection~\ref{f050404c})
for the class $\G^*$. In both cases we suppose that the group $G\in\G^*$ is given by a set
of generators. If $G\le\GL(n,R)$ for a certain $n\in\N$ and for a finite commutative ring $R$,
then in the case of LTP we set $V$ to be the standard free $R$-module of dimension $n$ on
which the group $\GL(n,R)$ acts, whereas for MT problem we set $F=\GL(n,R)$.

\lmml{f070405a}
Let $G\in\G^*$. Then given a derivation tree of $G$ the problems MT
and LTP can be solved in time polynomial in $L(G)$.
\elmm
\proof Let $T$ be a derivation tree of $G$. Then the labels of the leaves of $T$ are
the groups $G_1,\ldots,G_t\in\G_0$. Due to the choice of $\G_0$ the problems MT and LTP can be
solved for the group $G_i$ in time polynomial in $L(G_i)$ for $i=1,\ldots,t$.
Since $L(G)=L(T)^{O(1)}$, it suffices to verify that by means of the tree $T$ the problems
can be reduced in time $L(T)^{O(1)}$ to the corresponding problems for $G_1,\ldots,G_t$. For
this purpose
let us consider, for instance, the
reduction in the case of the primitive wreath product $G=H\wr\Gamma$ with $H\le\GL(n,R)$ and
$\Gamma=\sym(m)$ (other operations from $\O$ on groups are treated in a similar way). Then $G\le\GL(n^m,R)$
and since $T$ is given, we know the decomposition
$$
V=U\otimes \cdots \otimes U\quad (m\ \mbox{times})
$$
where $V$ and $U$ are the standard $R$-modules for groups $\GL(n^m,R)$ and $\GL(n,R)$
respectively. Any element $g\in G$ can be represented as the pair $(h,k)\in H^m\times\sym(m)$
such that
\qtnl{f210405b}
(u_1,\ldots,u_m)^g=(u_{i_1}^{h_{i_1}},\ldots,u_{i_m}^{h_{i_m}})
\eqtn
where $h=(h_1,\ldots,h_m)$ and $i_j=j^{k^{-1}}$ for $j=1,\ldots,m$. Now the permutation
$k$ can be efficiently computed from the elements of the form $(0_R,\ldots,1_R,\ldots,0_R)^g$
(with $1_R$ being the unique nonzero component in a certain place). So the element $h=gg_k^{-1}$
also can be
found efficiently where $g_k$ is the element of $\GL(V)=\GL(n^m,R)$ corresponding to~$k$
(this element acts on $V$ exactly by permuting coordinates according to~$k$). In particular,
this provides a polynomial time reduction of the MT problem for $G$ to the corresponding
problem for $H$.

Next, proceeding to the LTP problem
let $v\in u^G$ for some $u,v\in V$. Denote by $D$ the bipartite graph with parts
being the multisets $\{u_1,\ldots,u_m\}$ and $\{v_1,\ldots,v_m\}$ and the edges
being the pairs $(u_i,v_j)$ for which $v_i\in (u_j)^H$. Then from (\ref{f210405b}) it
follows that there is a one to one correspondence between the matchings
$\{(u_i,v_{j_i}):\ i=1,\ldots,m\}$ of the graph $D$ and the set $\{k\in\Gamma:\ v=u^g$ with
$g=(h,k)\in G$ for some $h\in H^m\}$. Since the problem of finding a matching of a bipartite
graph can be solved efficiently, we see that the LTP problem for $G$ is polynomial time
reducible to the corresponding problem for $H$.\bull

A natural way to apply our construction to the key agreement protocol is to choose a random
group $G\in\G^*$ of a prescribed size and then choose random subgroups $G_A$ and $G_B$ of
$G$ (see~(\ref{f010804c})). These groups can be specified by sets of generators constructed
as follows:
\vspace{2mm}

\noindent{\bf Step 1.} Let $S$ be the set of leaves of the derivation tree of the group~$G$.
For each $s\in S$ take random subsets $X_A(s)$ and $X_B(s)$ of the group $H_s$ associated
with~$s$.
\vspace{2mm}

\noindent{\bf Step 2.} Using the natural embedding $h\to g_h$ of $H_s$ into $G$ output
$X_A=\{g_x:\ x\in X_A(s),\ s\in S\}$ and $X_B=\{g_x:\ x\in X_B(s),\ s\in S\}$ as the generator
sets of $G_A$ and $G_B$ respectively.
\vspace{2mm}

\noindent
Thus, the constructing of the groups $G_A$ and $G_B$ is performed simultaneously with the
constructing the group $G$. (In fact, all we need is the embedding of each group assigned
to a leaf of the derivation tree of the group $G$ into $G$.) In
this way it is possible to control some properties of the groups, for instance, to avoid the
situation when $G_A$ centralizes $G_B$ (then the common key coincides with $1_G$ and so is not
secure).

Applying our construction to design homomorphic cryptosystems is more delicate. First of
all we define the set $M(G)$ for each group $G\le\GL(n,R)$ for some $n\in\N$ and some finite
commutative ring $R$ (note that this covers the case $G\in\G_0$ and also allows one to produce
homomorphisms in one more way: replacing $\G_0$ by a bigger subclass of $\G$).
Namely, any automorphism $\sigma\in\aut(R)$ induces a homomorphism
$$
f_\sigma:G\to G^\sigma,\ A\mapsto A^\sigma
$$
where the matrix $A^\sigma\in\GL(n,R)$ is obtained from the matrix $A\in\GL(n,R)$ by
entry-wise applying of~$\sigma$. To choose $\sigma$ we observe that $R=\oplus_{i\in I}R_i$
where each $R_i$ is a finite local commutative ring. Any automorphism of the residue field of the
ring $R_i$ can be lifted to the automorphism of this ring (statement (3) of
Proposition~\ref{f090505b}). In the representation of the Galois ring
as a quotient ring of a ring of polynomials this lifting can be done efficiently. Taking any
collection $\{\sigma_i\}_{i\in I}$
one can construct the automorphism $\sigma\in\aut(R)$ such that $\sigma|_{R_i}=\sigma_i$ for
all~$i$. The set of such automorphisms we denote by $\aut_0(R)$ (in the case of $R$ being a
field this group coincides with $\aut(R)$). Set
\qtnl{f090405d}
\M(G)=f_0\cup\{f_\sigma:\ \sigma\in\aut_0(R)\}
\eqtn
where $f_0$ is a trivial homomorphism taking any element of $G$ to the identity matrix
of $\GL(n,R)$. Then assuming that the ring $R$ is given explicitly, one can choose a
random element of $\M(G)$ in time polynomial in $L(G)$.

To provide the recursive step in constructing homomorphisms take $o\in\O_s$, $s\ge 1$.
Suppose first that $s=1$. Then $o$ is an unary operation, i.e. either it changes the
underlying ring $R$ of a group $G\le\GL(n,R)$, or $o$ is a conjugation.
Given a homomorphism $f:G\to H$
with $H\le\GL(n,R)$ we set $o(f)$ to be the composition $o\circ f$. Let now $s>1$
and $f_i:G_i\to H_i$ be a homomorphism with $G_i,H_i\in\G^*$ for $i=1,\ldots,s$. Then there
exists the natural canonical homomorphism
$$
o(f_1,\ldots,f_s):o(G_1,\ldots,G_s)\to o(H_1,\ldots,H_s)
$$
coinciding with $f_i$ on the group $G_i$ which in this case is a subgroup of the
group $o(G_1,\ldots,G_s)$. In any case, the resulting homomorphism is efficiently
computable (we recall that we represent a homomorphism by listing explicitly the
images of the generators). The above discussion shows that the following statement
holds.

\lmml{f080405a}
Let $f:G\to H$ be a homomorphism constructed in the above way where $G,H\in\G^*$. Then given
the derivation tree of $G$ one can find $f(g)$ for $g\in G$ in time polynomial in $L(G)$
and the size of $g$.\bull
\elmm

\sbsnt{Secure generation.}\label{f090405a}
Let us fix the classes $\G_0,\G,\G^*$, the set $\O$ of operations and the sets $\M(G)$ for
$G\in\G_0$ as in Subsection~\ref{f050405u}. Then due to Lemmas~\ref{f070405a}
and~\ref{f080405a} one can construct groups $G\in\G^*$ to realize both key agreement protocols
and homomorphic cryptosystems in which the group $G$ and the derivation tree $T$ of it play
the roles of public and secret keys, respectively. The security of such systems is based on
the difficulty of the following problem.
\vspace{2mm}

\prblm{Decomposition Problem}{Given a group $G\in\G^*$ find a derivation tree $T$ of $G$.}
\vspace{2mm}

This problem arises in connection with a computational version of the above mentioned
Aschbacher's theorem. A number of practical algorithms (without complexity bounds) for
Decomposition Problem are known
(see \cite{LG}) but in general this problem seems to be difficult. Indeed, suppose that
$R=\ZZ_m$ where $m=pq$ with $p$ and $q$ being two different primes. Denote by $G_p$
the cyclic matrix group of the order $p$ in $\GL(2,p)$ (see~(\ref{f210405a})). Similarly,
the group $G_q$ is defined. Then $G_p,G_q\in\G_0$ and
$$
G=G'_p\times G'_q\le\GL(2,R)
$$
where $G'_p$ and $G'_q$ are the images of the groups $G_p$ and $G_q$ with respect to the
natural embeddings  $\GL(2,p)$ and $\GL(2,q)$ into $\GL(2,R)$. Thus the group $G$ can be
constructed in two steps: construct the groups $G'_p$ and $G'_q$ (the operation of
changing the underlying ring), and set $G=G'_1\times G'_2$ (the operation of the direct
product). This implies that $G\in\G^*$. This shows that the integer factoring problem is a
special case of the Decomposition Problem.

Another strategy of Charlie could be to avoid solving the Decomposition Problem and to try solve
the problems like LTP, SCSP or CMT directly. To prevent such an attack one can choose
the leaves of the derivation tree of the group $G$ to be the groups of the exponential size
with respect to $L(G)$. Then from the construction it follows that these groups will
arise as the composition factors of $G$. However, for the groups with large composition
factors all the problems like LTP, SCSP or CMT seem to be difficult (see Subsections~\ref{f050404c}
and~\ref{f090405c}).

We mention one more attack of Charlie for the case of a homomorphic cryptosystem. Suppose
we construct in the above way the homomorphism $f:G\to H$ with $G,H\in\G^*$. We call the
homomorphism
{\it linear} if it induces the ring homomorphism $f':A(G)\to A(H)$ where $A(G)$ (resp. $A(H)$)
is the subring of the underlying full matrix ring generated by $G$ (resp. $H$). For a linear
homomorphism the corresponding homomorphic cryptosystem can be easily broken whenever
$G\le\GL(n,R)$ where $R=\ZZ_n$ for some $n\in\N$ or $R$ is a finite field (or, more generally,
a direct sum of Galois rings). Indeed, in
this case Charlie can find $f(g)$ for $g \in G$ as follows. Take random generators
$g_1,\ldots,g_s$ of the group $G$ and find a decomposition $g=\sum_{i=1}^sc_ig_i$ with $c_i\in R$
just involving linear algebra. Then $f(g)=\sum_{i=1}^sc_if(g_i)$ due to the linearity of~$f$.
To prevent this attack one can take some initial homomorphisms at the leaves of the
derivation tree to be elements of the group $\aut_0(R)$ (see (\ref{f090405d})). Then
the constructed homomorphism is not linear in general (e.g. if $g\in\GL(n,\F)$ with $\F$ being
a field, and $\sigma\in\aut(\F)$, then generally $(ag)^\sigma\ne ag^\sigma$).

We complete the subsection by the following statement summarizing the above discussion.

\thrm
Assuming that the problems LTP, SCSP, CMT for matrix groups over finite commutative rings,
as well as the Decomposition Problem are intractable, a secure two-party key agreement
protocol and homomorphic cryptosystem can be implemented for these groups.\bull
\ethrm

One of the consequences of this theorem is that by means of it one can construct encrypted
simulation of a boolean circuit of the logarithmic depth (the details can be found
in~\cite{GP}).

\section*{Final remarks}
One of the main problems in constructing homomorphic public-key cryptosystems consists
in finding appropriate trapdoor functions. However, in the natural presentations of
homomorphisms of algebraic structures the problem of breaking such a system is
reduced to some variants of the CMT problem. On the other hand, there is the following result
for matrix groups over finite fields.

\thrml{f101204a}{\rm \cite[Theorem~6.1]{KS}}
Given $K=\lg X\rg\le\GL(d,p^e)$ where $X\subset\GL(d,p^e)$, there is a Las Vegas algorithm that
given any $g\in\GL(d,p^e)$, decides whether $g\in K$, and if $g\in K$, then the algorithm
produces
a straight-line program with the input $X$, yielding g. The algorithm uses an oracle to
compute discrete logarithms in fields of characteristic $p$ with sizes up to
$p^{ed}$. In  case when all of those composition factors of Lie type in
characteristic $p$ are constructively recognizable with a Discrete Log oracle
\footnote{The current list of groups of Lie type recognizable with a Discrete Log oracle
is given in~\cite{KS}; this list includes the groups of series A, B, C, D.},
the running time is a polynomial in the input length $|X|d^2e\log p$, plus the
time required for polynomially many calls to the Discrete Log oracle.\bull
\ethrm

This theorem
shows that having an oracle for the Discrete Logarithm,
the membership problem can be solved in probabilistic polynomial time for matrix
groups over finite fields.
This means that at least for homomorphic public-key cryptosystems over such groups
there is a little hope to find a trapdoor function different from
functions the difficulty of inversion of which is based on the intractability of the Discrete
Logarithm.
However, only a little is known on the computational complexity of the membership
problem for matrix groups over rings. So constructions over such groups seems to be
more perspective from the point of view of algebraic (non-commutative) cryptography.

\end{document}